\input amstex

\input amstex
\input amsppt.sty
\magnification=\magstep1
\hsize=36truecc
 \vsize=23.5truecm
\baselineskip=14truept
 \NoBlackBoxes
\def\q{\quad}
\def\qq{\qquad}
\def\mod#1{\ (\text{\rm mod}\ #1)}

\def\t{\text}
\def\qtq#1{\q\t{#1}\q}
\def\mod#1{\ (\text{\rm mod}\ #1)}
\def\qtq#1{\q\t{#1}\q}
\def\f{\frac}
\def\e{\equiv}
\def\b{\binom}

\def\sls#1#2{(\f{#1}{#2})}
 \def\ls#1#2{\big(\f{#1}{#2}\big)}

\par\q
\let \pro=\proclaim
\let \endpro=\endproclaim
\topmatter
\title Some supercongruences modulo $p^2$
\endtitle
\author ZHI-Hong Sun\endauthor
\affil School of the Mathematical Sciences, Huaiyin Normal
University,
\\ Huaian, Jiangsu 223001, PR China
\\ Email: zhihongsun$\@$yahoo.com
\\ Homepage: http://www.hytc.edu.cn/xsjl/szh
\endaffil

 \nologo \NoRunningHeads

\abstract{Let $p>3$ be a prime, and let $m$ be an integer with
$p\nmid m$. In the paper we prove some supercongruences concerning
$$\align &\sum_{k=0}^{p-1}\f{\b{2k}k\b{3k}k}{54^k},\
\sum_{k=0}^{p-1}\f{\b{2k}k\b{4k}{2k}}{128^k},\
\sum_{k=0}^{p-1}\f{\b{3k}k\b{6k}{3k}}{432^k},
\\&\sum_{k=0}^{p-1}\f{\b{2k}k^2\b{3k}{k}}{m^k},
\sum_{k=0}^{p-1}\f{\b{2k}k^2\b{4k}{2k}}{m^k},\
\sum_{k=0}^{p-1}\f{\b{2k}k\b{3k}{k}\b{6k}{3k}}{m^k}\mod
{p^2}.\endalign$$ Thus we solve some conjectures of Zhi-Wei Sun and
the author.
\par\q
 }
 \endabstract
  \footnote"" {The author is
supported by the Natural Sciences Foundation of China (grant No.
10971078).}
\endtopmatter
\document
\subheading{1. Introduction}
\par For positive integers
$a,b$ and $n$, if $n=ax^2+by^2$ for some integers $x$ and $y$, we
briefly say that $n=ax^2+by^2$. Let $p>3$ be a prime. In 2003,
Rodriguez-Villegas[RV] posed some conjectures on supercongruences
modulo $p^2$. Three of his conjectures are equivalent to
$$\align
&\sum_{k=0}^{p-1}\f{\b{2k}k^2\b{3k}k}{108^k}\e
 \cases 4A^2-2p\mod{p^2}&\t{if $p=A^2+3B^2\e 1\mod 3$,}
\\0\mod{p^2}&\t{if $p\e 2\mod 3$,}\endcases\tag 1.1
\\&\sum_{k=0}^{p-1}\f{\b{2k}k^2\b{4k}{2k}}{256^k}\e
 \cases 4c^2-2p\mod{p^2}&\t{if $p=c^2+2d^2\e 1,3\mod 8$,}
\\0\mod{p^2}&\t{if $p\e 5,7\mod 8$,}\endcases\tag 1.2
\\&\sum_{k=0}^{p-1}\f{\b{2k}k\b{3k}k\b{6k}{3k}}{1728^k}\e
\cases \sls p3(4a^2-2p)\mod{p^2}&\t{if $p=a^2+b^2\e 1\mod 4$ and
$2\nmid a$,}
\\0\mod{p^2}&\t{if $p\e 3\mod 4$,}\endcases\tag 1.3
\endalign$$
where  $\sls am$ is the Jacobi symbol. The above conjectures have
been solved by Mortenson[M] and Zhi-Wei Sun[Su2].
\par Let $\Bbb Z$ be the set of
integers, and let $[x]$ be the greatest integer function. For a
prime $p$ let $\Bbb Z_p$ be the set of rational numbers whose
denominator is coprime to $p$. Recently the author's brother Zhi-Wei
Sun posed many conjectures ([Su1]) involving
$$\align &\sum_{k=0}^{p-1}\f{\b{2k}k\b{3k}k}{54^k},\
\sum_{k=0}^{p-1}\f{\b{2k}k\b{4k}{2k}}{128^k},\
\sum_{k=0}^{p-1}\f{\b{3k}k\b{6k}{3k}}{432^k},
\\&\sum_{k=0}^{p-1}\f{\b{2k}k^2\b{3k}{k}}{m^k},
\sum_{k=0}^{p-1}\f{\b{2k}k^2\b{4k}{2k}}{m^k},\
\sum_{k=0}^{p-1}\f{\b{2k}k\b{3k}{k}\b{6k}{3k}}{m^k}\mod
{p^2},\endalign$$ where $p>3$ is a prime and $m\in\Bbb Z$ with
$p\nmid m$. For example, Zhi-Wei Sun conjectured ([Su1, Conjectures
A8 and A9]) that for any prime $p>3$,
$$\align
&\sum_{k=0}^{p-1}\f{\b{2k}k^2\b{3k}k}{(-192)^k}\e\cases
0\mod{p^2}&\t{if $p\e 2\mod 3$,}
\\L^2-2p\mod {p^2}&\t{if $p\e 1\mod 3$ and so $4p=L^2+27M^2$,}\endcases
\tag 1.4
\\&\sum_{k=0}^{p-1}\f{(6k)!}{(-96)^{3k}(3k)!k!^3}
\e\cases 0\mod {p^2}&\t{if $\sls p{19}=-1$,}
\\\sls{-6}p(x^2-2p)\mod {p^2}&\t{if $\sls p{19}=1$ and so $4p=x^2+19y^2$.}
\endcases\tag 1.5\endalign$$ In [S3], the author proved (1.4) and
(1.5) modulo $p$.
\par Let $p$ be an odd prime and let $x$ be a variable.
In the paper we establish the following general congruences:
$$\align &\sum_{k=0}^{p-1}\b{2k}k^2\b{3k}k(x(1-27x))^k\e \Big(
\sum_{k=0}^{p-1}\b{2k}k\b{3k}kx^k\Big)^2\mod {p^2},
\\&\sum_{k=0}^{p-1}\b{2k}k^2\b{4k}{2k}(x(1-64x))^k\e \Big(
\sum_{k=0}^{p-1}\b{2k}k\b{4k}{2k}x^k\Big)^2\mod {p^2},
\\&\sum_{k=0}^{p-1}\b {2k}k\b{3k}k\b{6k}{3k}(x(1-432x))^k\e \Big(
\sum_{k=0}^{p-1}\b {3k}k\b{6k}{3k}x^k\Big)^2\mod {p^2}.
\endalign$$
As an application, using the work in [S2,S3] we prove many
congruences modulo $p^2$. For example, (1.4) is true for $p\e 2\mod
3$ and (1.5) is true when $\sls p{19}=-1$.

 \subheading{2. Congruences for $\sum_{k=0}^{p-1}
 \f{\b{2k}k^2\b{3k}k}{m^k}$ and
  $\sum_{k=0}^{p-1}\f{\b{2k}k^2\b{4k}{2k}}{m^k}\mod{p^2}$}
\pro{Lemma 2.1} Let $m$ be a nonnegative integer. Then
$$\sum_{k=0}^m\b {2k}k^2\b{3k}k\b k{m-k}(-27)^{m-k}
=\sum_{k=0}^m\b {2k}k\b{3k}k\b{2(m-k)}{m-k}\b{3(m-k)}{m-k}.$$
\endpro We prove the lemma by using WZ method and Mathematica.
Clearly the result is true for $m=0,1$. Since both sides satisfy the same recurrence relation
$$81(m+1)(3m+2)(3m+4)S(m)
    -3(2m+3)(9m^2+27m+22)S(m+1) + (m+2)^3 S(m+2) = 0,$$
    we see that the lemma is true.
The  proof certificate for the left hand side is
$$ - \frac{729 k^2(m+2)(m-2k)(m-2k+1)}{(m-k+1)(m-k+2)},
$$ and the proof certificate for the right hand side is $$
\frac{9
k^2(3m-3k+1)(3m-3k+2)(9m^2-9mk+30m-14k+24)}{(m-k+1)^2(m-k+2)^2}. $$
\pro{Theorem 2.1} Let $p$ be an odd prime and let  $x$ be a
variable. Then
$$\sum_{k=0}^{p-1}\b{2k}k^2\b{3k}k(x(1-27x))^k\e \Big(
\sum_{k=0}^{p-1}\b{2k}k\b{3k}kx^k\Big)^2\mod {p^2}.$$
\endpro
Proof. It is clear that
$$\align &\sum_{k=0}^{p-1}\b{2k}k^2\b{3k}k(x(1-27x))^k
\\&=\sum_{k=0}^{p-1}\b{2k}k^2\b{3k}kx^k\sum_{r=0}^k\b kr(-27x)^r
\\&=\sum_{m=0}^{2(p-1)}x^m\sum_{k=0}^{min\{m,p-1\}}\b{2k}k^2\b{3k}k\b
k{m-k}(-27)^{m-k}.\endalign$$
 Suppose $p\le m\le 2p-2$ and $0\le
k\le p-1$. If $k>\f p2$, then $p\mid \b{2k}k$ and so $p^2\mid
\b{2k}k^2$. If $k<\f p2$, then $m-k\ge p-k>k$ and so $\b k{m-k}=0$.
Thus, from the above and Lemma 2.1 we deduce
$$\align &\sum_{k=0}^{p-1}\b{2k}k^2\b{3k}k(x(1-27x))^k
\\&\e \sum_{m=0}^{p-1}x^m\sum_{k=0}^m\b {2k}k^2\b{3k}k
\b k{m-k}(-27)^{m-k}
\\&=\sum_{m=0}^{p-1}x^m\sum_{k=0}^m\b {2k}k\b{3k}k
\b{2(m-k)}{m-k}\b{3(m-k)}{m-k}
\\&=\sum_{k=0}^{p-1}\b{2k}k\b{3k}kx^k\sum_{m=k}^{p-1}
\b{2(m-k)}{m-k}\b{3(m-k)}{m-k}x^{m-k}
\\&=\sum_{k=0}^{p-1}\b{2k}k\b{3k}kx^k\sum_{r=0}^{p-1-k}\b{2r}r\b{3r}rx^r
\\&=\sum_{k=0}^{p-1}\b{2k}k\b{3k}kx^k\Big(\sum_{r=0}^{p-1}\b{2r}r\b{3r}rx^r
-\sum_{r=p-k}^{p-1}\b{2r}r\b{3r}rx^r\Big)
\\&=\Big(\sum_{k=0}^{p-1}\b{2k}k\b{3k}kx^k\Big)^2
-\sum_{k=0}^{p-1}\b{2k}k\b{3k}kx^k
\sum_{r=p-k}^{p-1}\b{2r}r\b{3r}rx^r
 \mod{p^2}.\endalign$$
If $\f{2p}3\le k\le p-1$, then $\b{2k}k\b{3k}k=\f{(3k)!}{k!^3}\e
0\mod {p^2}$. If $0\le k \le \f p3$ and $p-k\le r\le p-1$, then
$\f{2p}3\le r\le p-1$ and so $\b{2r}r\b{3r}r=\f{(3r)!}{r!^3}\e 0\mod
{p^2}$. If $\f p3<k<\f{2p}3$ and $p-k\le r\le p-1$, then $r\ge
p-k>\f p3$, $\b{2k}k\b{3k}k=\f{(3k)!}{k!^3}\e 0\mod p$ and
$\b{2r}r\b{3r}r=\f{(3r)!}{r!^3}\e 0\mod p$. Hence, for $0\le k\le
p-1$ and $p-k\le r\le p-1$ we have $p^2\mid
\b{2k}k\b{3k}k\b{2r}r\b{3r}r$ and so
$$\sum_{k=0}^{p-1}\b{2k}k\b{3k}kx^k
\sum_{r=p-k}^{p-1}\b{2r}r\b{3r}rx^r
 \e 0\mod{p^2}.$$
Therefore the result follows.

\pro{Corollary 2.1} Let $p>3$ be a prime and $m\in\Bbb Z_p$ with
$m\not\e 0\mod p$. Then
$$\sum_{k=0}^{p-1}\f{\b{2k}k^2\b{3k}k}{m^k}
\e\Big(\sum_{k=0}^{p-1}\b{2k}k\b{3k}k\Big(
\f{1-\sqrt{1-108/m}}{54}\Big)^k\Big)^2 \mod {p^2}.$$
\endpro
Proof. Taking $x=\f{1-\sqrt{1-108/m}}{54}$ in Theorem 2.1 we deduce
the result.
\pro{Corollary 2.2} Let $p>3$ be a prime and $m\in\Bbb
Z_p$ with $m\not\e 0\mod p$. Then
$$\sum_{k=0}^{[p/3]}\f{\b{2k}k^2\b{3k}k}{m^k}\e 0\mod p\qtq{implies}
\sum_{k=0}^{p-1}\f{\b{2k}k^2\b{3k}k}{m^k}\e 0\mod {p^2}.$$
\endpro
Proof. Clearly $\b{2k}k\b{3k}k=\f{(3k)!}{k!^3}\e 0\mod p$ for $\f
p3<k<p$. Suppose $\sum_{k=0}^{[p/3]}\f{\b{2k}k^2\b{3k}k}{m^k}\e
0\mod p$. Then
$$\sum_{k=0}^{p-1}\f{\b{2k}k^2\b{3k}k}{m^k}\e
\sum_{k=0}^{[p/3]}\f{\b{2k}k^2\b{3k}k}{m^k}\e 0\mod p.$$ Using
Corollary 2.1 we see that
$$\sum_{k=0}^{p-1}\b{2k}k\b{3k}k\Big(
\f{1-\sqrt{1-108/m}}{54}\Big)^k\e 0\mod p$$ and so the result
follows from Corollary 2.1.

 \pro{Theorem 2.2} Let
$p\e 1\mod 3$ be a prime and so $p=A^2+3B^2$ with $A\e 1\mod 3$.
Then
$$\sum_{k=0}^{p-1}\f{\b{2k}k\b{3k}k}{54^k}
\e 2A-\f p{2A}\mod{p^2}.$$
\endpro
Proof. Clearly $\b{2k}k\b{3k}k=\f{(3k)!}{k!^3}\e 0\mod p$ for $\f
p3<k<p$. Using [S1, Theorem 2.5] we have
$$\sum_{k=0}^{p-1}\f{\b{2k}k\b{3k}k}{54^k}=
\sum_{k=0}^{p-1}\f{(3k)!}{54^k\cdot k!^3}\e
 \sum_{k=0}^{[p/3]}\f{(3k)!}{54^k\cdot k!^3}\e 2A\mod p.$$
 Set $\sum_{k=0}^{p-1}\f{\b{2k}k\b{3k}k}{54^k}=2A+qp$. Then
 $$\Big(\sum_{k=0}^{p-1}\f{\b{2k}k\b{3k}k}{54^k}\Big)^2=(2A+qp)^2\e
 4A^2+4Aqp\mod{p^2}.$$ Taking $m=108$ in Corollary 2.1 we get
 $$\sum_{k=0}^{p-1}\f{\b{2k}k^2\b{3k}k}{108^k}
 \e \Big(\sum_{k=0}^{p-1}\f{\b{2k}k\b{3k}k}{54^k}\Big)^2\mod{p^2}.$$
  Thus, by (1.1) and the above we have
 $$4A^2-2p\e \sum_{k=0}^{p-1}\f{\b{2k}k^2\b{3k}k}{108^k}\e
 \Big(\sum_{k=0}^{p-1}\f{\b{2k}k\b{3k}k}{54^k}\Big)^2\e 4A^2+4Aqp\mod{p^2}$$ and hence
 $q\e -\f 1{2A}\mod p$. So the theorem is proved.
 \newline{\bf Remark 2.1} Theorem 2.2 was conjectured by the author
 in [S1]. When $p$ is a prime of the form $6k+5$, it was conjectured in [S1] that
 $\sum_{k=0}^{p-1}\f{\b{2k}k\b{3k}k}{54^k}\e 0\mod{p^2}$. This was
 recently confirmed by Zhi-Wei Sun[Su3].

 \pro{Theorem 2.3} Let $p\e 5\mod 6$ be a prime.
 Then
 $$\sum_{k=0}^{p-1}\f{\b{2k}k^2\b{3k}k}{(-192)^k}\e 0\mod{p^2}.$$
 \endpro
Proof. Since $\b{2k}k\b{3k}k=\f{(3k)!}{k!^3}\e 0\mod p$ for $\f
p3<k<p$, from [S3, Theorem 4.3] we know that
$$\sum_{k=0}^{p-1}\f{\b{2k}k\b{3k}k}{(-216)^k}
=\sum_{k=0}^{p-1}\f{(3k)!}{(-216)^k\cdot
 k!^3}\e \sum_{k=0}^{[p/3]}\f{(3k)!}{(-216)^k\cdot
 k!^3}\e 0\mod p.$$
 Thus, taking $x=-\f 1{216}$ in Theorem 2.1 we obtain the result.

 \pro{Lemma 2.2} Let $m$ be a nonnegative integer. Then
$$\sum_{k=0}^m\b {2k}k^2\b{4k}{2k}\b k{m-k}(-64)^{m-k}
=\sum_{k=0}^m\b {2k}k\b{4k}{2k}\b{2(m-k)}{m-k}\b{4(m-k)}{2(m-k)}.$$
\endpro
We prove the lemma by using WZ method and Mathematica. Clearly the
result is true for $m=0,1$. Since both sides satisfy the same
recurrence relation $$ \align &1024(m+1)(2m+1)(2m+3)S(m)
    -8(2m+3)(8m^2+24m+19)S(m+1)\\&\qq+(m+2)^3 S(m+2) = 0,\endalign$$ we see that
    Lemma 2.2 is true.
 The  proof certificate for the left hand side is
$$ - \frac{4096 k^2(m+2)(m-2k)(m-2k+1)}{(m-k+1)(m-k+2)},$$
 and the proof certificate for the right hand side is
$$ \frac{16
k^2(4m-4k+1)(4m-4k+3)(16m^2-16mk+55m-26k+46)}{(m-k+1)^2(m-k+2)^2}.$$
\pro{Theorem 2.4} Let $p$ be an odd prime and let  $x$ be a
variable. Then
$$\sum_{k=0}^{p-1}\b{2k}k^2\b{4k}{2k}(x(1-64x))^k\e \Big(
\sum_{k=0}^{p-1}\b{2k}k\b{4k}{2k}x^k\Big)^2\mod {p^2}.$$
\endpro
Proof. It is clear that
$$\align &\sum_{k=0}^{p-1}\b{2k}k^2\b{4k}{2k}(x(1-64x))^k
\\&=\sum_{k=0}^{p-1}\b{2k}k^2\b{4k}{2k}x^k\sum_{r=0}^k\b kr(-64x)^r
\\&=\sum_{m=0}^{2(p-1)}x^m\sum_{k=0}^{min\{m,p-1\}}\b{2k}k^2\b{4k}{2k}\b
k{m-k}(-64)^{m-k}.\endalign$$
 Suppose $p\le m\le 2p-2$ and $0\le
k\le p-1$. If $k>\f p2$, then $p\mid \b{2k}k$ and so $p^2\mid
\b{2k}k^2$. If $k<\f p2$, then $m-k\ge p-k>k$ and so $\b k{m-k}=0$.
Thus, from the above and Lemma 2.2 we deduce
$$\align &\sum_{k=0}^{p-1}\b{2k}k^2\b{4k}{2k}(x(1-64x))^k
\\&\e \sum_{m=0}^{p-1}x^m\sum_{k=0}^m\b {2k}k^2\b{4k}{2k}
\b k{m-k}(-64)^{m-k}
\\&=\sum_{m=0}^{p-1}x^m\sum_{k=0}^m\b {2k}k\b{4k}{2k}
\b{2(m-k)}{m-k}\b{4(m-k)}{2(m-k)}
\\&=\sum_{k=0}^{p-1}\b{2k}k\b{4k}{2k}x^k\sum_{m=k}^{p-1}
\b{2(m-k)}{m-k}\b{4(m-k)}{2(m-k)}x^{m-k}
\\&=\sum_{k=0}^{p-1}\b{2k}k\b{4k}{2k}x^k
\sum_{r=0}^{p-1-k}\b{2r}r\b{4r}{2r}x^r
\\&=\sum_{k=0}^{p-1}\b{2k}k\b{4k}{2k}x^k\Big(\sum_{r=0}^{p-1}\b{2r}r
\b{4r}{2r}x^r -\sum_{r=p-k}^{p-1}\b{2r}r\b{4r}{2r}x^r\Big)
\\&=\Big(\sum_{k=0}^{p-1}\b{2k}k\b{4k}{2k}x^k\Big)^2
-\sum_{k=0}^{p-1}\b{2k}k\b{4k}{2k}x^k
\sum_{r=p-k}^{p-1}\b{2r}r\b{4r}{2r}x^r
 \mod{p^2}.\endalign$$
Now suppose $0\le k\le p-1$ and $p-k\le r\le p-1$. If $k\ge
\f{3p}4$, then $p^2\nmid (2k)!$, $p^3\mid (4k)!$ and so $
\b{2k}k\b{4k}{2k}=\f{(4k)!}{(2k)!k!^2}\e 0\mod{p^2}$. If $k < \f
p4$, then $r\ge p-k\ge \f{3p}4$ and so $
\b{2r}r\b{4r}{2r}=\f{(4r)!}{(2r)!r!^2}\e 0\mod{p^2}$. If $\f p4<k<\f
p2$, then $r\ge p-k>\f p2$, $p\nmid (2k)!$, $p\mid (4k)!$, $p\mid
\b{2r}r$ and $\b{2k}k\b{4k}{2k}=\f{(4k)!}{(2k)!k!^2}\e 0\mod p$. If
$\f p2<k<\f{3p}4$, then $r\ge p-k>\f p4$, $p\mid \b{2k}k$ and
$\b{2r}r\b{4r}{2r}=\f{(4r)!}{(2r)!r!^2}\e 0\mod p$. Hence we always
have $\b{2k}k\b{4k}{2k}\b{2r}r\b{4r}{2r}\e 0\mod{p^2}$ and so
$$\sum_{k=0}^{p-1}\b{2k}k\b{4k}{2k}x^k
\sum_{r=p-k}^{p-1}\b{2r}r\b{4r}{2r}x^r
 \e 0\mod{p^2}.$$
Now combining all the above we obtain the result. \pro{Corollary
2.3} Let $p>3$ be a prime and $m\in\Bbb Z_p$ with $m\not\e 0\mod p$.
Then
$$\sum_{k=0}^{p-1}\f{\b{2k}k^2\b{4k}{2k}}{m^k}
\e\Big(\sum_{k=0}^{p-1}\b{2k}k\b{4k}{2k}\Big(
\f{1-\sqrt{1-256/m}}{128}\Big)^k\Big)^2 \mod {p^2}.$$
\endpro
Proof. Taking $x=\f{1-\sqrt{1-256/m}}{128}$ in Theorem 2.4 we deduce
the result. \pro{Corollary 2.4} Let $p>3$ be a prime and $m\in\Bbb
Z_p$ with $m\not\e 0\mod p$. Then
$$\sum_{k=0}^{[p/4]}\f{\b{2k}k^2\b{4k}{2k}}{m^k}\e 0\mod p\qtq{implies}
\sum_{k=0}^{p-1}\f{\b{2k}k^2\b{4k}{2k}}{m^k}\e 0\mod {p^2}.$$
\endpro
Proof. For $\f p4<k<p$ we see that
$\b{2k}k^2\b{4k}{2k}=\f{(4k)!}{k!^4}\e 0\mod p$. Suppose
$\sum_{k=0}^{[p/4]}\f{\b{2k}k^2\b{4k}{2k}}{m^k}\e 0\mod p$. Then
$$\sum_{k=0}^{p-1}\f{\b{2k}k^2\b{4k}{2k}}{m^k}\e
\sum_{k=0}^{[p/4]}\f{\b{2k}k^2\b{4k}{2k}}{m^k}\e 0\mod p.$$ Using
Corollary 2.3 we see that
$$\sum_{k=0}^{p-1}\b{2k}k\b{4k}{2k}\Big(
\f{1-\sqrt{1-256/m}}{128}\Big)^k\e 0\mod p$$ and so the result
follows from Corollary 2.3.

 \pro{Theorem 2.5}
Let $p\e 1,3\mod 8$ be a prime and $p=c^2+2d^2$ with $c,d\in\Bbb Z$
and $c\e 1\mod 4$. Then
$$\sum_{k=0}^{p-1}\f{\b{2k}k\b{4k}{2k}}{128^k}
\e (-1)^{[\f p8]+\f{p-1}2}\Big(2c-\f p{2c}\Big)\mod{p^2}.$$
\endpro
Proof. From the proof of Theorem 2.4 we know that $p\mid
\b{2k}k\b{4k}{2k}$ for $\f p4<k<p$. By [S2, Theorem 2.1] we have
$$\sum_{k=0}^{p-1}\f{\b{2k}k\b{4k}{2k}}{128^k}\e
\sum_{k=0}^{[p/4]}\f{\b{2k}k\b{4k}{2k}}{128^k}
 \e (-1)^{[\f p8]+\f{p-1}2}2c\mod p.$$
  Set $\sum_{k=0}^{p-1}\f{\b{2k}k\b{4k}{2k}}{128^k}
=(-1)^{[\f p8]+\f{p-1}2}2c+qp$. Then
 $$\Big(\sum_{k=0}^{p-1}\f{\b{2k}k\b{4k}{2k}}{128^k}\Big)^2
 =((-1)^{[\f p8]+\f{p-1}2}2c+qp)^2\e
 4c^2+(-1)^{[\f p8]+\f{p-1}2}4cqp\mod{p^2}.$$ Taking $m=256$
  in Corollary 2.3 we get
 $$\sum_{k=0}^{p-1}\f{\b{2k}k^2\b{4k}{2k}}{256^k}
 \e \Big(\sum_{k=0}^{p-1}\f{\b{2k}k\b{4k}{2k}}{128^k}\Big)^2\mod{p^2}.$$
  Thus, by (1.2) and the above,
 $$4c^2-2p\e  \sum_{k=0}^{p-1}\f{\b{2k}k^2\b{4k}{2k}}{256^k}
 \e\Big(\sum_{k=0}^{p-1}\f{\b{2k}k\b{4k}{2k}}{128^k}\Big)^2
 \e 4c^2+(-1)^{[\f p8]+\f{p-1}2}4cqp\mod{p^2}$$ and hence
 $q\e -(-1)^{[\f p8]+\f{p-1}2}\f 1{2c}\mod p$. So the theorem is proved.
 \par\q
 \newline{\bf Remark 2.2} Theorem 2.5 is a conjecture of Zhi-Wei
 Sun ([Su1, Conjecture A49]). In [Su3], Zhi-Wei Sun showed that
$\sum_{k=0}^{p-1}\f{\b{2k}k\b{4k}{2k}}{128^k}\e 0\mod{p^2}$ for
primes $p\e 5,7\mod 8$.
 \pro{Theorem 2.6 ([S1, Conjecture 2.1])} Let $p>3$ be a prime of
the form $4k+3$. Then
$$\sum_{k=0}^{p-1}\f{\b{2k}k^2\b{4k}{2k}}{648^k}\e 0\mod{p^2}.$$
\endpro
Proof. Since $p\mid \b{2k}k\b{4k}{2k}$ for $p>k>\f p4$, from [S2,
Theorem 2.4] we know that
$$\sum_{k=0}^{p-1}\f{\b{2k}k\b{4k}{2k}}{72^k}
\e \sum_{k=0}^{[p/4]}\f{\b{2k}k\b{4k}{2k}}{72^k}\e 0\mod p.$$
 Thus, taking $x=\f 1{72}$ in Theorem 2.4 and applying the above
 we obtain the result.

 \pro{Theorem 2.7 ([S1, Conjecture 2.2])} Let $p$ be a prime of the form $6k+5$. Then
$$\sum_{k=0}^{p-1}\f{\b{2k}k^2\b{4k}{2k}}{(-144)^k}\e 0\mod{p^2}.$$
\endpro
Proof. Since $p\mid \b{2k}k\b{4k}{2k}$ for $p>k>\f p4$, from [S2,
Theorem 2.5] we know that
$$\sum_{k=0}^{p-1}\f{\b{2k}k\b{4k}{2k}}{48^k}
\e \sum_{k=0}^{[p/4]}\f{\b{2k}k\b{4k}{2k}}{48^k}\e 0\mod p.$$
 Thus, taking $x=\f 1{48}$ in Theorem 2.4 and applying the above
 we obtain the result.

  \pro{Theorem 2.8 ([S1, Conjecture 2.3])} Let $p>3$ be a prime such that $p\e 3,5,6
  \mod 7$. Then
$$\sum_{k=0}^{p-1}\f{\b{2k}k^2\b{4k}{2k}}{(-3969)^k}\e 0\mod{p^2}.$$
\endpro
Proof. Since $p\mid \b{2k}k\b{4k}{2k}$ for $p>k>\f p4$, from [S2,
Theorem 2.6] we know that
$$\sum_{k=0}^{p-1}\f{\b{2k}k\b{4k}{2k}}{63^k}
\e \sum_{k=0}^{[p/4]}\f{\b{2k}k\b{4k}{2k}}{63^k}\e 0\mod p.$$
 Thus, taking $x=\f 1{63}$ in Theorem 2.4 and applying the above
 we obtain the result.
\subheading{3. Congruences for
$\sum_{k=0}^{p-1}\f{\b{2k}k\b{3k}k\b{6k}{3k}}{m^k}\mod{p^2}$}

  \pro{Lemma 3.1} Let $m$ be a nonnegative
integer. Then
$$\sum_{k=0}^m\b {2k}k\b{3k}k\b{6k}{3k}\b k{m-k}(-432)^{m-k}
=\sum_{k=0}^m\b {3k}k\b{6k}{3k}\b{3(m-k)}{m-k}\b{6(m-k)}{3(m-k)}.$$
\endpro
\par We prove the lemma by using WZ method and Mathematica. Clearly the result is true
 for $m=0,1$. Since both sides satisfy the same recurrence relation $$
\align &20736(m+1)(3m+1)(3m+5)S(m)-24(2m+3)(18m^2+54m+41)S(m+1)
    \\&\qq+(m+2)^3 S(m+2) = 0,\endalign$$ we see that Lemma 3.1 is true.
The  proof certificate for the left hand side is $$
 - \frac{186624 k^2(m+2)(m-2k)(m-2k+1)}{(m-k+1)(m-k+2)},$$
and the proof certificate for the right hand side is $$ \frac{144
k^2(6m-6k+1)(6m-6k+5)(36m^2-36mk+129m-62k+114)}{(m-k+1)^2(m-k+2)^2}.$$
\par For given prime $p$ and integer $n$, if $p^{\alpha}\mid n$ but
$p^{\alpha+1}\nmid n$, we say that $p^{\alpha}\ \Vert\ n$.
\pro{Lemma 3.2} Let $p$ be an odd prime, $k,r\in\{0,1,\ldots, p-1\}$
and $k+r\ge p$. Then
$$\b{3k}k\b{6k}{3k}\b{3r}r\b{6r}{3r}\e 0\mod{p^2}.$$
\endpro
Proof. If $k>\f{5p}6$, then $p^5\mid (6k)!$, $p\ \Vert\ (2k)!$,
$p^2\ \Vert\ (3k)!$ and so
$\b{3k}k\b{6k}{3k}=\f{(6k)!}{k!(2k)!(3k)!}\e 0\mod{p^2}$. If
$\f{2p}3\le k<\f{5p}6$, then $2p\le  3k<3p$, $4p\le 6k<5p$, $p^4\
\Vert\  (6k)!$, $p^2\ \Vert\ (3k)!$, $p\ \Vert\ (2k)!$ and so
$\b{3k}k\b{6k}{3k}=\f{(6k)!}{k!(2k)!(3k)!}\e 0\mod p$. If $\f
p2<k<\f{2p}3$, then $p<3k<2p$, $3p<6k<4p$, $p^3\mid (6k)!$, $p\
\Vert\ (2k)!$, $p\ \Vert\ (3k)!$ and so
$\b{3k}k\b{6k}{3k}=\f{(6k)!}{k!(2k)!(3k)!}\e 0\mod p$. If $\f p3\le
k<\f p2$, then $2k<p$, $p\le 3k<2p$, $6k\ge 2p$, $p^2\mid (6k)!$,
$p\nmid (2k)!$, $p\ \Vert\ (3k)!$ and so
$\b{3k}k\b{6k}{3k}=\f{(6k)!}{k!(2k)!(3k)!}\e 0\mod p$. If $\f
p6<k<\f p3$, then $3k<p$, $6k>p$ and so
$\b{3k}k\b{6k}{3k}=\f{(6k)!}{k!(2k)!(3k)!}\e 0\mod p$.
\par From the above we see that $p\mid \b{3k}k\b{6k}{3k}$ for $k>\f
p 6$. Therefore, if $k>\f p6$ and $r>\f p6$, then
$\b{3k}k\b{6k}{3k}\b{3r}r\b{6r}{3r}\e 0\mod{p^2}.$ If $r<\f p6$,
then $k\ge p-r>\f{5p}6$ and so $p^2\mid \b{3k}k\b{6k}{3k}$ by the
above. If $k<\f p6$, then $r\ge p-k>\f{5p}6$ and so $p^2\mid
\b{3r}r\b{6r}{3r}$ by the above.
\par Now putting all the above together we prove the lemma.

 \pro{Theorem 3.1} Let $p$ be an odd prime and let  $x$ be a
variable. Then
$$\sum_{k=0}^{p-1}\b {2k}k\b{3k}k\b{6k}{3k}(x(1-432x))^k\e \Big(
\sum_{k=0}^{p-1}\b {3k}k\b{6k}{3k}x^k\Big)^2\mod {p^2}.$$
\endpro
Proof. It is clear that
$$\align &\sum_{k=0}^{p-1}\b {2k}k\b{3k}k\b{6k}{3k}(x(1-432x))^k
\\&=\sum_{k=0}^{p-1}\b {2k}k\b{3k}k\b{6k}{3k}x^k\sum_{r=0}^k\b kr(-432x)^r
\\&=\sum_{m=0}^{2(p-1)}x^m\sum_{k=0}^{min\{m,p-1\}}\b {2k}k\b{3k}k\b{6k}{3k}
\b k{m-k}(-432)^{m-k}.\endalign$$
 Suppose $p\le m\le 2p-2$ and $0\le
k\le p-1$. If $k\ge \f {2p}3$, then $2p\le 3k<3p$, $6k\ge 4p$,
$p^3\nmid (3k)!$, $p^4\mid (6k)!$ and so $\b
{2k}k\b{3k}k\b{6k}{3k}=\f{(6k)!}{(3k)!k!^3}\e 0\mod {p^2}$. If $\f
p2<k<\f{2p}3$, then $3k<2p$, $6k>3p$, $p^2\nmid (3k)!$ and $p^3\mid
(6k)!$ and so $\b {2k}k\b{3k}k\b{6k}{3k}=\f{(6k)!}{(3k)!k!^3}\e
0\mod {p^2}$.  If $k<\f p2$, then $m-k\ge p-k>k$ and so $\b
k{m-k}=0$. Thus, from the above and Lemma 3.1 we deduce
$$\align &\sum_{k=0}^{p-1}\b {2k}k\b{3k}k\b{6k}{3k}(x(1-432x))^k
\\&\e \sum_{m=0}^{p-1}x^m\sum_{k=0}^m\b {2k}k\b{3k}k\b{6k}{3k}
\b k{m-k}(-432)^{m-k}
\\&=\sum_{m=0}^{p-1}x^m\sum_{k=0}^m\b {3k}k\b{6k}{3k}
\b{3(m-k)}{m-k}\b{6(m-k)}{3(m-k)}
\\&=\sum_{k=0}^{p-1}\b{3k}k\b{6k}{3k}x^k\sum_{m=k}^{p-1}
\b{3(m-k)}{m-k}\b{6(m-k)}{3(m-k)}x^{m-k}
\\&=\sum_{k=0}^{p-1}\b{3k}k\b{6k}{3k}x^k
\sum_{r=0}^{p-1-k}\b{3r}r\b{6r}{3r}x^r
\\&=\sum_{k=0}^{p-1}\b{3k}k\b{6k}{3k}x^k\Big(\sum_{r=0}^{p-1}\b{3r}r
\b{6r}{3r}x^r -\sum_{r=p-k}^{p-1}\b{3r}r\b{6r}{3r}x^r\Big)
\\&=\Big(\sum_{k=0}^{p-1}\b{3k}k\b{6k}{3k}x^k\Big)^2
-\sum_{k=0}^{p-1}\b{3k}k\b{6k}{3k}x^k
\sum_{r=p-k}^{p-1}\b{3r}r\b{6r}{3r}x^r
 \mod{p^2}.\endalign$$
 By Lemma 3.2, we have $p^2\mid \b{3k}k\b{6k}{3k}\b{3r}r\b{6r}{3r}$
 for $0\le k\le p-1$ and $p-k\le r\le p-1$.
Thus
$$\sum_{k=0}^{p-1}\b{3k}k\b{6k}{3k}x^k
\sum_{r=p-k}^{p-1}\b{3r}r\b{6r}{3r}x^r
  \e 0\mod{p^2}.$$
 Now combining all the above we obtain the result.

 \pro{Corollary 3.1} Let $p>3$ be a prime and $m\in\Bbb Z_p$ with
$m\not\e 0\mod p$. Then
$$\sum_{k=0}^{p-1}\f{\b{2k}k\b{3k}k\b{6k}{3k}}{m^k}
\e\Big(\sum_{k=0}^{p-1}\b{3k}k\b{6k}{3k}\Big(
\f{1-\sqrt{1-1728/m}}{864}\Big)^k\Big)^2 \mod {p^2}.$$
\endpro
Proof. Taking $x=\f{1-\sqrt{1-1728/m}}{864}$ in Theorem 3.1 we
deduce the result. \pro{Corollary 3.2} Let $p>3$ be a prime and
$m\in\Bbb Z_p$ with $m\not\e 0\mod p$. Then
$$\sum_{k=0}^{[p/6]}\f{\b{2k}k\b{3k}k\b{6k}{3k}}{m^k}\e 0\mod p\qtq{implies}
\sum_{k=0}^{p-1}\f{\b{2k}k\b{3k}k\b{6k}{3k}}{m^k}\e 0\mod {p^2}.$$
\endpro
Proof. From the proof of Lemma 3.2 we know that $p\mid
\b{3k}k\b{6k}{3k}$ for $p>k>\f p6$. Suppose
$\sum_{k=0}^{[p/6]}\f{\b{2k}k\b{3k}k\b{6k}{3k}}{m^k}\e 0\mod p$.
Then
$$\sum_{k=0}^{p-1}\f{\b{2k}k\b{3k}k\b{6k}{3k}}{m^k}\e
\sum_{k=0}^{[p/6]}\f{\b{2k}k\b{3k}k\b{6k}{3k}}{m^k}\e 0\mod p.$$
Using Corollary 3.1 we see that
$$\sum_{k=0}^{p-1}\b{3k}k\b{6k}{3k}\Big(
\f{1-\sqrt{1-1728/m}}{864}\Big)^k\e 0\mod p$$ and so the result
follows. \pro{Theorem 3.2} Let $p\e 1\mod 4$ be a prime and
$p=a^2+b^2$ with $a,b\in\Bbb Z$ and $4\mid a-1$. Then
$$\sum_{k=0}^{p-1}\f{\b{6k}{3k}\b{3k}k}{864^k}\e \cases
2a-\f p{2a}\mod {p^2}&\t{if $p\e 1\mod {12}$ and $3\nmid a$,} \\
-2a+\f p{2a}\mod {p^2}&\t{if $p\e 1\mod {12}$ and $3\mid a$,}
\\2b-\f p{2b} \mod {p^2}&\t{if $p\e 5\mod{12}$ and $3\mid a-b$.}
\endcases$$
\endpro
Proof. By the proof of Lemma 3.2, we have $p\mid \b{3k}k\b{6k}{3k}$
for $p>k>\f p6$. Set $$r=\cases a&\t{if $p\e 1\mod{12}$ and $3\nmid
a$,}
 \\-a&\t{if $p\e 1\mod{12}$ and $3\mid a$,}
 \\b&\t{if $p\e 5\mod{12}$ and $3\mid a-b$}.\endcases$$
 Using [S3, Theorem 2.1] we have
$$\sum_{k=0}^{p-1}\f{\b{3k}k\b{6k}{3k}}{864^k}\e
\sum_{k=0}^{[p/6]}\f{\b{3k}k\b{6k}{3k}}{864^k}
 \e 2r\mod p.$$
  Set $\sum_{k=0}^{p-1}\f{\b{3k}k\b{6k}{3k}}{864^k}
=2r+qp$. Then
 $$\Big(\sum_{k=0}^{p-1}\f{\b{3k}k\b{6k}{3k}}{864^k}\Big)^2
 =(2r+qp)^2\e
 4r^2+4rqp\mod{p^2}.$$ Taking $m=1728$
  in Corollary 3.1 we get
 $$\sum_{k=0}^{p-1}\f{\b{2k}k\b{3k}k\b{6k}{3k}}{1728^k}
 \e \Big(\sum_{k=0}^{p-1}\f{\b{3k}k\b{6k}{3k}}{864^k}\Big)^2\mod{p^2}.$$
  Thus, by (1.3) and the above,
 $$\ls p3(4a^2-2p)\e  \sum_{k=0}^{p-1}\f{\b{2k}k\b{3k}k\b{6k}{3k}}{1728^k}
 \e\Big(\sum_{k=0}^{p-1}\f{\b{3k}k\b{6k}{3k}}{864^k}\Big)^2
 \e 4r^2+4rqp\mod{p^2}$$ and hence
 $q\e -\f 1{2r}\mod p$. So the theorem is proved.
\par\q
\newline{\bf Remark 3.1} Theorem 3.2 is a conjecture of Zhi-Wei Sun
([Su1, Conjecture A44]). In [Su3], Zhi-Wei Sun showed that
$\sum_{k=0}^{p-1}\f{\b{3k}k\b{6k}{3k}}{864^k}\e 0\mod{p^2}$ for any
prime $p=4n+3>3$.

\pro{Theorem 3.3 ([S1, Conjectures 2.8 and 2.9])} Let $p>7$ be a
prime such that $p\e 3,5,6\mod 7$. Then
$$\sum_{k=0}^{p-1}\f{\b{2k}k\b{3k}k\b{6k}{3k}}{(-15)^{3k}}
\e \sum_{k=0}^{p-1}\f{\b{2k}k\b{3k}k\b{6k}{3k}}{255^{3k}}\e
0\mod{p^2}.$$
\endpro
Proof. This is immediate from [S3, Theorems 3.9 and 3.10] and
Corollary 3.2.

\pro{Theorem 3.4 ([Su1, Conjecture A26])} Let $p$ be an odd prime
with $p\e 2,6,7,8,$ $10\mod{11}$. Then
$$\sum_{k=0}^{p-1}\f{\b{2k}k\b{3k}k\b{6k}{3k}}{(-32)^{3k}}\e
0\mod{p^2}.$$
\endpro
Proof. This is immediate from [S3, Theorem 3.4] and Corollary 3.2.

\pro{Theorem 3.5 ([Su1, Conjecture A9])} Let $p>3$ be a prime with
$\sls p{19}=-1$. Then
$$\sum_{k=0}^{p-1}\f{\b{2k}k\b{3k}k\b{6k}{3k}}{(-96)^{3k}}\e
0\mod{p^2}.$$
\endpro
Proof. This is immediate from [S3, Theorem 3.5] and Corollary 3.2.

\pro{Theorem 3.6 ([Su1, Conjecture A10])} Let $p>5$ be a prime with
$\sls p{43}=-1$. Then
$$\sum_{k=0}^{p-1}\f{\b{2k}k\b{3k}k\b{6k}{3k}}{(-960)^{3k}}\e
0\mod{p^2}.$$
\endpro
Proof. This is immediate from [S3, Theorem 3.6] and Corollary 3.2.

\pro{Theorem 3.7 ([Su1, Conjecture A11])} Let $p>5$ be a prime with
$p\not=11$ and $\sls p{67}=-1$. Then
$$\sum_{k=0}^{p-1}\f{\b{2k}k\b{3k}k\b{6k}{3k}}{(-5280)^{3k}}\e
0\mod{p^2}.$$
\endpro
Proof. This is immediate from [S3, Theorem 3.7] and Corollary 3.2.

\pro{Theorem 3.8 ([Su1, Conjecture A12])} Let $p\not=2,3,5,23,29$ be
a prime with $\sls p{163}=-1$. Then
$$\sum_{k=0}^{p-1}\f{\b{2k}k\b{3k}k\b{6k}{3k}}{(-640320)^{3k}}\e
0\mod{p^2}.$$
\endpro
Proof. This is immediate from [S3, Theorem 3.8] and Corollary 3.2.

\pro{Theorem 3.9 ([S1, Conjecture 2.4])} Let $p\not=3,11$ be a prime
such that $p\e 3\mod 4$. Then
$$\sum_{k=0}^{p-1}\f{\b{2k}k\b{3k}k\b{6k}{3k}}{66^{3k}}
\e 0\mod{p^2}.$$
\endpro
Proof. This is immediate from [S3, Theorem 3.11] and Corollary 3.2.

\pro{Theorem 3.10 ([S1, Conjecture 2.5])} Let $p>5$ be a prime such
that $p\e 5,7\mod 8$. Then
$$\sum_{k=0}^{p-1}\f{\b{2k}k\b{3k}k\b{6k}{3k}}{20^{3k}}
\e 0\mod{p^2}.$$
\endpro
Proof. This is immediate from [S3, Theorem 3.12] and Corollary 3.2.

\pro{Theorem 3.11 ([S1, Conjecture 2.6])} Let $p>5$ be a prime such
that $p\e 2\mod 3$. Then
$$\sum_{k=0}^{p-1}\f{\b{2k}k\b{3k}k\b{6k}{3k}}{54000^k}
\e 0\mod{p^2}.$$
\endpro
Proof. This is immediate from [S3, Theorem 3.13] and Corollary 3.2.
\par\q
\newline{\bf Acknowledgements.} The author is indebted to Prof. Qing-Hu
Hou for his help in proving Lemmas 2.1, 2.2 and 3.1 by using WZ
method.

\enddocument

\Refs \widestnumber\key {BEW}

 \ref\key M\by  E. Mortenson\paper
Supercongruences for truncated $\ _{n+1}F_n$ hypergeometric series
with applications to certain weight three newforms\jour Proc. Amer.
Math. Soc.\vol 133(2005)\pages 321-330.\endref

 \ref\key RV\by F. Rodriguez-Villegas
\paper  Hypergeometric families of Calabi-Yau manifolds. Calabi-Yau
Varieties and Mirror Symmetry (Yui, Noriko (ed.) et al., Toronto,
ON, 2001), 223-231, Fields Inst. Commun., 38, Amer. Math. Soc.,
Providence, RI, 2003\endref


 \ref\key S1\by Z.H. Sun\paper
Congruences concerning Legendre polynomials\jour Proc. Amer. Math.
Soc. (2010), doi: 10.1090/S0002-9939-2010-10566-X, arXiv:1012.3833.
http://arxiv.org/abs/1012.3833
\endref

\ref\key S2\by Z.H. Sun\paper Congruences concerning Legendre
polynomials II\jour Amer. J. Math., submitted\finalinfo
\newline arXiv:1012.3898. http://arxiv.org/abs/1012.3898\endref
\ref\key S3\by Z.H. Sun\paper Congruences concerning Legendre
polynomials III\jour Acta Arith., submitted\finalinfo
\newline arXiv:1012.4234. http://arxiv.org/abs/1012.4234\endref

 \ref \key Su1\by Z.W. Sun\paper
Open conjectures on congruences, arXiv:0911.5665v53.
http://arxiv.org/abs/0911.5665\endref

\ref \key Su2\by Z.W. Sun\paper On sums involving products of three
binomial coefficients, arXiv:1012.3141. http://\newline
arxiv.org/abs/1012.3141\endref

\ref \key Su3\by Z.W. Sun\paper Super congruences and elliptic
curves over $\Bbb F_p$, arXiv:1011.6676. http://arxiv.org\newline
/abs/1011.6676\endref


\endRefs
\bye